\documentclass[11pt,a4paper]{amsart}
\usepackage[T1]{fontenc}
\usepackage{amssymb}
\usepackage[margin=25mm]{geometry}
\usepackage{microtype}
\usepackage{cite}
\usepackage[colorlinks=true,linkcolor=blue,citecolor=blue,urlcolor=blue]{hyperref}
\numberwithin{equation}{section}
\newtheorem{theorem}{Theorem}[section]
\newtheorem{lemma}{Lemma}[section]
\allowdisplaybreaks[2]
\title[Mixed modular and Dirichlet moments]{Lower bounds for mixed moments of modular
and Dirichlet $L$-functions}
\author{Yu Hao Qiu and Nanxiang Wang*}
\address{Department of Mathematics, China University of Mining and Technology,
Beijing 100083, People's Republic of China}
\email{YuHao.Qiu.China@hotmail.com, bqt2600701002@student.cumtb.edu.cn}
\address{School of Mathematics and Statistics, Qingdao University, Qingdao, Shandong 266071, People's Republic of China}
\email{wangnanxiang@qdu.edu.cn}
\date{}
\begin{document}
\begin{abstract}
Let $f$ be a fixed holomorphic Hecke cusp form of level one.
For every fixed positive integer $k$, we prove a lower bound of order
$q(\log q)^{2k^2}$ for the $2k$-th moment of
$L(1/2,f\otimes\chi)L(1/2,\chi)$ over primitive characters modulo
a sufficiently large prime $q$.
\end{abstract}
\maketitle

\section{Introduction}\label{hm:sec:introduction}

Moments of central values of $L$-functions describe their size and
distribution and provide a natural setting for nonvanishing questions.

Let $f$ be a fixed holomorphic Hecke cusp form of even weight $\kappa$
for $\mathrm{SL}_2(\mathbb Z)$, with Fourier expansion
\[
 f(z)=\sum_{n\ge1}\lambda_f(n)n^{(\kappa-1)/2}e(nz),
\]
with normalized $\lambda_f(1)=1.$ The coefficients $\lambda_f(n)$ are real. Deligne's bound
\cite{Deligne1974} gives
\begin{equation}\label{hm:eq:deligne}
 |\lambda_f(n)|\le d(n)\qquad(n\ge1),
\end{equation}
where $d(n)$ is the divisor function. For a Dirichlet character $\chi$
modulo $q$, write
\[
 \begin{aligned}
 L(s,f\otimes\chi)
 &=\sum_{n\ge1}\frac{\lambda_f(n)\chi(n)}{n^s}
 =\prod_p\left(1-\frac{\lambda_f(p)\chi(p)}{p^s}
                  +\frac{\chi(p)^2}{p^{2s}}\right)^{-1},\\
 L(s,\chi)&=\sum_{n\ge1}\frac{\chi(n)}{n^s},
 \qquad \Re s>1.
 \end{aligned}
\]
The first $L$-function extends to an entire function. If $\chi$ is
primitive, its completed function 
\[
\Lambda(s,f\otimes\chi)
 =\left(\frac q{2\pi}\right)^s
   \Gamma\!\left(s+\frac{\kappa-1}{2}\right)L(s,f\otimes\chi)
\]
satisfies the following functional equation: 
\[
 \Lambda(s,f\otimes\chi)
 =\epsilon(f\otimes\chi)\Lambda(1-s,f\otimes\overline\chi). 
\]
For the even characters ,
$\epsilon(f\otimes\chi)=i^\kappa\epsilon_\chi^2$, where
$\epsilon_\chi=\tau(\chi)/\sqrt q$ and $\tau(\chi)$ is the Gauss sum;
see \cite[Section~2.2]{BFKMM2023}.

The method of Rudnick and Soundararajan
\cite{RudnickSoundararajan2005,RudnickSoundararajan2006}
provides lower bounds for moments by introducing suitable short
Dirichlet polynomials. Chen and Jiang \cite{ChenJiang2023} applied
this approach to positive even moments of twists of holomorphic and
Maass forms for admissible moduli. For prime and prime-power moduli,
their estimates have the conjectured logarithmic order. Gao, He,
and Wu \cite{GaoHeWu2022} obtained sharp lower bounds for real
moments for prime powers with sufficiently large prime moduli, which using
the lower-bound principle
of Heap and Soundararajan \cite{HeapSoundararajan2022} and ideas from
the upper-bound method of Radziwi\l\l\ and Soundararajan
\cite{RadziwillSoundararajan2015}.

For products of two distinct twisted modular $L$ functions, Chen and
He \cite{ChenHe2021} proved lower bounds for even  moments. More generally, Chen, Li, and Wang
\cite{ChenLiWang2025} treated real moment parameters greater than
$1$ in their setting. Consider Conjugate mixed moments, it have been used to study simultaneous
nonvanishing; see Das and Khan \cite{DasKhan2015} and Zacharias
\cite{Zacharias2019}. 

Recently, asymptotic formulas and upper bounds for moments of products
of $L$-functions have been obtained by Chen and Lin \cite{ChenLin2026},
Hagen \cite{Hagen2026},  \cite{MQW2025},
\cite{TW2025}.
A related problem concerns averages over Dirichlet characters of the product
\[
L\!\left(\frac{1}{2}, f\otimes\chi\right)
L\!\left(\frac{1}{2}, \chi\right).
\]

Motivated by \cite{ChenJiang2023}, This paper aims to establish the following results.

\begin{theorem}\label{hm:thm:main}
For every fixed positive integer $k$ and every sufficiently large prime $q$,
\begin{equation}\label{thm:main-theorem}
 \sum_{\chi\bmod q}^{*}
 \left|L\!\left(\tfrac12,f\otimes\chi\right)
        L\!\left(\tfrac12,\chi\right)\right|^{2k}
 \gg_{f,k} q(\log q)^{2k^2}.
\end{equation}
The star denotes summation over primitive characters.
\end{theorem}

\noindent\textbf{Notations} Throughout, $q$ tends to infinity through odd primes, $\chi_0$ denotes
the principal character, and
\[
 \mathcal X_q^+
 =\{\chi\bmod q:\chi\ne\chi_0,\ \chi(-1)=1\}.
\]
All characters in this set are primitive. The symbol $\varepsilon$
denotes an arbitrarily small positive number and may change from one
occurrence to the next. We use $\ll$ and $\gg$ to denote Vinogradov's well-known notations. 
We write $A\asymp B$ when $A\ll B$ and
$B\ll A$. Bold letters denote vectors. The implied constants may
depend on the fixed form $f$, the fixed integer $k$, and any displayed
auxiliary parameters. We use $\overline a$ for the inverse of a nonzero
residue $a$ modulo $q$. $e(z)=e^{2\pi iz}.$

\section{Preliminaries}\label{hm:sec-preliminaries}
The following partition of unity permits a dyadic decomposition into
smooth sums; see \cite[Lemma~2.27]{BFKMM2023}.

\begin{lemma}\label{hm:lemma-dyadic}
There exists a nonnegative smooth function $W$, supported on
$[1/2,2]$, such that $W^{(j)}(x)\ll_j1$ for every integer $j\ge0$ and
\[
\sum_{v\ge0}W\!\left(\frac{x}{2^v}\right)=1
\qquad (x\ge1).
\]
\end{lemma}

\begin{lemma}[Unbalanced approximate functional equation]\cite[Lemma 2.2]{ChenJiang2023}\label{lem:unbalanced_afe}
For any real number \(Z \geq 1\), we have
\begin{equation}
L(1/2, f \otimes \chi) = \sum_{n \geq 1} \frac{\lambda(n)\chi(n)}{\sqrt{n}} U\left(\frac{n}{qZ}\right) + \epsilon(f \otimes \chi) \sum_{n \geq 1} \frac{\lambda(n)\bar{\chi}(n)}{\sqrt{n}} U\left(\frac{nZ}{q}\right),
\end{equation}
where
\begin{equation}
U(x) := \frac{1}{2\pi i} \int_{(2)} \frac{L_\infty(1/2 + s, f \otimes \chi)}{L_\infty(1/2, f \otimes \chi)} x^{-s} \frac{ds}{s}, 
\end{equation}
\[
L_{\infty}(s, f \otimes \chi)=(2\pi)^{-s}\Gamma((\kappa-1)/2+s)
\]
satisfies \(U^{(j)}(x) \ll_{j,A} (1 + x)^{-A}\) for all \(A,\ j \geq 0\), \(\epsilon(f \otimes \chi)=i^\kappa\epsilon_\chi^2\), where
$\epsilon_\chi=\tau(\chi)/\sqrt q$ and $\tau(\chi)$ is the Gauss sum.

\end{lemma}
\begin{lemma}\label{hm:lem:weighted}
Let $1\le h,\ \ell\le q^{1/16}$, q is a prime, and write
\[
g=(h,\ell),\qquad a=h/g,\qquad b=\ell/g.
\]
Define
\[
F(\chi)=L(1/2,f\otimes\chi)L(1/2,\overline\chi).
\]
Then, uniformly in $h$ and $\ell$,
\begin{equation}\label{hm:eq-weighted-first}
\sum_{\chi\in\mathcal X_q^+}
F(\chi)\chi(h)\overline{\chi(\ell)}
=\frac{q-1}{2\sqrt{ab}}
\sum_{\substack{r\ge1\\(r,q)=1}}\frac{\lambda_f(br)}r
+O_{f,\varepsilon}\!\left(q^{11/12+\varepsilon}\right),
\end{equation}
where the series $\sum\limits_{\substack{r\ge1\\(r,q)=1}}\frac{\lambda_f(br)}r$ is a convergent series. 
\end{lemma}

\begin{proof}
Put $Q=q^{3/2}$ and $H=\max(h,\ell)$. Since $\chi$ is even,
$\epsilon_\chi\epsilon_{\overline\chi}=1$. Thus, the root number of
the product that defines $F(\chi)$ is
$i^\kappa\epsilon_\chi^2\epsilon_{\overline\chi}
=i^\kappa\epsilon_\chi$.
Let $D_\chi(s)=L(s,f\otimes\chi)L(s,\overline{\chi})$, $\gamma(s)=L_\infty(s,f\otimes\chi)\pi^{-s/2}\Gamma\left(\frac{s}{2}\right)$,  $\Lambda_\chi(s)=q^s\gamma(s)D_\chi(s)$. For even primitive characters, define $\Lambda_D(s,\chi)=q^{s/2}\pi^{-s/2}\Gamma\left(\frac{s}{2}\right)L(s,\chi)$, it satisfies
$$\Lambda_D(s,\overline{\chi})=\epsilon_{\overline{\chi}}\Lambda_D(1-s,\chi).$$
Since
\[
 \Lambda(s,f\otimes\chi)
 =\epsilon(f\otimes\chi)\Lambda(1-s,f\otimes\overline\chi), 
\]
and $\chi$ is an even character, we have 
\begin{equation}\label{1}
    \Lambda_\chi(s)=i^\kappa\epsilon_\chi^2\epsilon_{\overline{\chi}}\Lambda_{\overline{\chi}}(1-s)=i^\kappa\epsilon_\chi\Lambda_{\overline{\chi}}(1-s).
\end{equation}
Define $U(x):=\frac{1}{2\pi i}\int_{(2)}\frac{\gamma\left(\frac{1}{2}+s\right)}{\gamma\left(\frac{1}{2}\right)}x^{-s}\frac{ds}{s}$. Consider the integral
\[I_\chi:=\frac{1}{2\pi i}\int_{(2)}Q^s\frac{\gamma\left(\frac{1}{2}+s\right)}{\gamma\left(\frac{1}{2}\right)}D_\chi\left(\frac{1}{2}+s\right)\frac{ds}{s}.\] 
On the vertical line $\mathfrak{R}s=2$, $L(s,f\otimes\chi)$ and $L(s,\overline{\chi})$ are convergent absolutely, hence 
\[D_\chi\left(\frac{1}{2}+s\right)=\sum_{n,m\ge1}\frac{\lambda_f(n)\chi(n)\overline{\chi(m)}}{(nm)^{1/2+s}}.\]
Commute the sum and the integral in $I_\chi$, we have 
\begin{equation}\label{2}
    I_\chi=\sum_{n,m\ge1}\frac{\lambda_f(n)\chi(n)\overline{\chi(m)}}{\sqrt{nm}}U\left(\frac{nm}{Q}\right).
\end{equation}
For $I_\chi$, note that between the line $\mathfrak{R}s=2$ and $\mathfrak{R}s=-2$, there is one pole $s=0$. Let 
\[J_\chi:=\frac{1}{2\pi i}\int_{(-2)}Q^s\frac{\gamma\left(\frac{1}{2}+s\right)}{\gamma\left(\frac{1}{2}\right)}D_\chi\left(\frac{1}{2}+s\right)\frac{ds}{s}\]
By the Residue Theorem, we have 
\[I_\chi=D_\chi\left(\frac{1}{2}\right)+J_\chi=F(\chi)+J_\chi.\]
By (\ref{1}), we have
\[Q^s\gamma\left(\frac{1}{2}+s\right)D_\chi\left(\frac{1}{2}+s\right)=i^\kappa Q^{-s}\gamma\left(\frac{1}{2}-s\right)D_{\overline{\chi}}\left(\frac{1}{2}-s\right).\] 
By changing the variable as $s\to-s$, we have 
\[J_\chi=i^\kappa\epsilon_\chi\frac{1}{2\pi i}\int_{(-2)}Q^s\frac{\gamma\left(\frac{1}{2}-s\right)}{\gamma\left(\frac{1}{2}\right)}D_{\overline{\chi}}\left(\frac{1}{2}-s\right)\frac{ds}{s}=-i^\kappa\epsilon_\chi I_{\overline{\chi}}.\]
Thus, $F(\chi)=I(\chi)+i^\kappa\epsilon_\chi I_{\overline{\chi}}.$ By (\ref{2}), we have 
\begin{equation}\label{hm:eq-product-afe}
\begin{aligned}
F(\chi)
=\sum_{n,m\ge1}
\frac{\lambda_f(n)\chi(n)\overline{\chi(m)}}{\sqrt{nm}}
U\!\left(\frac{nm}{Q}\right)+i^\kappa\epsilon_\chi\sum_{n,m\ge1}
\frac{\lambda_f(n)\overline{\chi(n)}\chi(m)}{\sqrt{nm}}
U\!\left(\frac{nm}{Q}\right),
\end{aligned}
\end{equation}
where the product weight is
\[
U(x)=\frac{1}{2\pi i}\int_{(2)}\frac{\gamma\left(\frac{1}{2}+s\right)}{\gamma\left(\frac{1}{2}\right)}x^{-s}\frac{ds}{s}=\frac{1}{2\pi i}\int_{(2)}(2\pi)^{-s}\pi^{-s/2}
\frac{\Gamma(\kappa/2+s)}{\Gamma(\kappa/2)}
\frac{\Gamma(1/4+s/2)}{\Gamma(1/4)}
x^{-s}\,\frac{ds}{s}.
\]
By Lemma \ref{lem:unbalanced_afe}, for every integer $j\ge0$ and every $B>0$, we have $x^jU^{(j)}(x)\ll_{f,j,B}(1+x)^{-B}$.

After multiplying \eqref{hm:eq-product-afe} by
$\chi(h)\overline{\chi(\ell)}$ and summing over
$\mathcal X_q^+$, denote the direct and dual contributions by
$D_q(h,\ell)$ and $R_q(h,\ell)$, respectively. Applying the following formula
\[
\sum_{\substack{\chi\bmod q\\\chi(-1)=1}}\chi(u)\overline{\chi(v)}=\frac{1}{2}\sum_{\chi\bmod q}(1+\chi(-1))\chi(u)\overline{\chi(v)}
\]
we derive that
\begin{equation}\label{hm:eq-direct-split}
D_q(h,\ell)=D_{\mathrm{diag}}(h, \ell)+D_{\mathrm{off}}(h, \ell).
\end{equation}
where
\[
D_{\mathrm{diag}}(h,l)=
\frac{q-1}{2}
\sum_{\substack{m,n\ge1 \\ (mn,q)=1 \\ nh \equiv ml \pmod q\\nh=ml}}
\frac{\lambda_f(n)}{\sqrt{mn}} U\!\left(\frac{mn}{Q}\right),
\]
and
\[
\begin{aligned}
D_{\mathrm{off}}(h,l)
={}& \frac{q-1}{2}
\sum_{\substack{m,n\ge 1\\ (mn,q)=1\\ nh\equiv ml\pmod q\\ nh\neq ml}}
\frac{\lambda_f(n)}{\sqrt{mn}}\, U\!\left(\frac{mn}{Q}\right) \\[2mm]
&+\frac{q-1}{2}
\sum_{\substack{m,n\ge 1\\ (mn,q)=1\\ nh\equiv -ml\pmod q}}
\frac{\lambda_f(n)}{\sqrt{mn}}\, U\!\left(\frac{mn}{Q}\right) \\[2mm]
&-\sum_{\substack{m,n\ge 1\\ (mn,q)=1}}
\frac{\lambda_f(n)}{\sqrt{mn}}\, U\!\left(\frac{mn}{Q}\right).
\end{aligned}
\]
Next, we need to deal with $D_{\mathrm{diag}}(h, l)$. Indeed, the condition $nh=m\ell$ is equivalent to
$n=br$, $m=ar$. So we can rewrite $D_{\mathrm{diag}}(h,l)$ into
\[
    D_{\mathrm{diag}}(h,l)=\frac{q-1}{2}\sum_{\substack{r\geq 1\\(r, q)=1}}\frac{\lambda_f(br)}{r} U\!\left(\frac{abr^2}{Q}\right).
\]
Next, we need to prove $\sum\limits_{\substack{r\ge1\\(r,q)=1}}\frac{\lambda_f(br)}r$ converges. Firstly, using the well-known result
\[
\sum_{d|(m, n)} \mu(d) = 
\begin{cases} 
1, & (m,n) = 1; \\ 
0, & (m,n) \geq 2.
\end{cases}
\]
we can rewrite this series into 
\[
\sum\limits_{\substack{r\ge1\\(r,q)=1}}\frac{\lambda_f(br)}r=\sum_{d\mid q}\mu(d)\sum_{\substack{r\geq 1\\d\mid r}}\frac{\lambda_f(br)}{r}=\sum_{d\mid q}\mu(d)\sum_{k\geq 1}\frac{\lambda_f(bkd)}{kd}=\sum_{d\mid q}\frac{\mu(d)}{d}\sum_{k\geq 1}\frac{\lambda_f(bkd)}{k}.
\]
Then according to the relation
\[
\lambda_f(km)=\sum_{l\mid(m,k)}\mu(l)\lambda_f(m/l)\lambda_f(k/l),
\]
let $bd=lt_1, k=lt_2,$ then we can rewrite
\[
\sum_{k\leqslant X}\frac{\lambda_f(bkd)}{k}=\sum_{l\leqslant X}\frac{\mu(l)\lambda_f(t_1)}{l}\sum_{t_2\leqslant X/l}\frac{\lambda_f(t_2)}{t_2}
\]
Next, using Abel's  summation, the bound $|\lambda_f(n)|\le d(n)\ll n^\epsilon$ and the best-known result due to Rankin\cite{R89}: 
\[
\sum_{n\leqslant x}\lambda_f(n)\ll_f x^{1/3}(\log x)^{-(1-8/(3\pi))},
\]
we have
\[\sum_{k\leqslant X}\frac{\lambda_f(bkd)}{k}\ll_f(bd)^{\varepsilon}\]

Then according to Cauchy's criterion of series, we complete the proof.
Then substituting the 
representation of $U$ into $D_{\mathrm{diag}}$ and moving the contour
to $\Re s=-1/4$ crosses only the pole at $s=0$.
The new integral is therefore
\[
\ll_{f,\varepsilon}
\frac q{\sqrt{ab}}\left(\frac Q{ab}\right)^{-1/4}q^\varepsilon
\ll_{f,\varepsilon}q^{5/8+\varepsilon}.
\]
It follows that
\begin{equation}\label{hm:eq-diagonal-main}
D_{\mathrm{diag}}(h, l)=
\frac{q-1}{2\sqrt{ab}}
\sum_{\substack{r\ge1\\(r,q)=1}}\frac{\lambda_f(br)}r
+O_{f,\varepsilon}\!\left(q^{5/8+\varepsilon}\right).
\end{equation}

To estimate $D_{\mathrm{off}}$, restrict to a dyadic block
$n\asymp N$, $m\asymp M$, and denote $P=NM$. The number of solutions
to either congruence $nh\equiv\pm m\ell\pmod q$ is
$O(P/q+\min(N,M))$. If such a solution does not lie on the ordinary
diagonal, then $\max(N,M)\gg q/H$. Thus, using
$|\lambda_f(n)|\le d(n)$, the contribution of such a block is bounded by
\[
\begin{aligned}
&\ll_{f,\varepsilon,B}
P^\varepsilon\left(\sqrt P+
q\sqrt{\frac{\min(N,M)}{\max(N,M)}}\right)(1+P/Q)^{-B}\\
&\ll_{f,\varepsilon,B}H P^{1/2+\varepsilon}(1+P/Q)^{-B}.
\end{aligned}
\]
Taking $B$ sufficiently large and summing the dyadic blocks gives
\begin{equation}\label{hm:eq-direct-error}
D_{\mathrm{off}}(h,l)\ll_{f,\varepsilon}Hq^{3/4+\varepsilon}
\ll_{f,\varepsilon}q^{13/16+\varepsilon}.
\end{equation}

We next estimate the dual contribution. Expanding the Gauss sum and
applying character's  orthogonality gives
\[
\sum_{\chi\in\mathcal X_q^+}
\epsilon_\chi\chi(u)\overline{\chi(v)}
=\frac{q-1}{2\sqrt q}
\left[e\!\left(\frac{v\overline u}{q}\right)
+e\!\left(-\frac{v\overline u}{q}\right)\right]
+\frac1{\sqrt q}.
\]
The positive correction arises from subtracting the principal Gauss
sum $\tau(\chi_0)=-1$. Taking $u=hm$, $v=\ell n$, and
$c\equiv\ell\overline h\pmod q$, we obtain
\[
\begin{aligned}
R_q(h,\ell)
={}&\frac{i^\kappa(q-1)}{2\sqrt q}
\sum_{\substack{n,m\ge1\\(nm,q)=1}}
\frac{\lambda_f(n)}{\sqrt{nm}}U\!\left(\frac{nm}{Q}\right)
\left[e\!\left(\frac{cn\overline m}{q}\right)
+e\!\left(-\frac{cn\overline m}{q}\right)\right]\\
&+\frac{i^\kappa}{\sqrt q}
\sum_{\substack{n,m\ge1\\(nm,q)=1}}
\frac{\lambda_f(n)}{\sqrt{nm}}U\!\left(\frac{nm}{Q}\right)
=R_{1q}(h,l)+R_{2q}(h,l).
\end{aligned}
\]
By \eqref{hm:eq:deligne}, we obtain
\begin{equation}\label{hm:eq-dual-correction}
R_{2q}(h,l)
\ll_f q^{-1/2}Q^{1/2}(\log(2Q))^2
\ll_{f,\varepsilon}q^{1/4+\varepsilon}.
\end{equation}
To estimate $R_{1q}(h,l)$, it suffices to bound
\[
E(N,M;c)=\frac{\sqrt q}{\sqrt{NM}}
\left|\sum_{\substack{n,m\ge1\\(nm,q)=1}}
\lambda_f(n)e\!\left(\frac{cn\overline m}{q}\right)
W_{N,M}\!\left(\frac nN,\frac mM\right)\right|.
\]
where 
\[
W_{N,M}(n, m)=U\left(\frac{nm}{Q}\right)\frac{1}{\sqrt{nm}}.
\]
The weight function $W_{N,M}$ is supported on $[1/2,2]^2$; we can derive that
\[
\frac{\partial^{i+j}}{\partial x^i \partial y^j} W_{N,M}(x,y) \ll_{f,i,j,B} \left(1 + \frac{P}{Q}\right)^{-B}
\]
where $P=NM$.
Then use the following result
\[
\sum_{n\le Y}\lambda_f(n)e(\alpha n)
\ll_f Y^{1/2}\log(2Y);
\]
see \cite[equation~(1.2)]{Godber2013}.
 Summing over \(n\) first, and then summing over \(m\), we obtain
\begin{equation}\label{hm:eq-wilton-block}
E(N,M;c)
\ll_{f,B}q^{1/2}M^{1/2}\log(2N)(1+P/Q)^{-B}.
\end{equation}

  Alternatively fix $n$, let $W_n(y)=W_{N,M}(n/N,y)$ and define
$\widehat W_n(\xi)=\int_{\mathbb R}W_n(y)e(-y\xi)\,dy$.
By Poisson summation, we obtain
\[
\sum_{\substack{m\in\mathbb Z\\q\nmid m}}
W_n(m/M)e\!\left(\frac{cn\overline m}{q}\right)
=\frac Mq\sum_{r\in\mathbb Z}S(r,cn;q)\widehat W_n(rM/q),
\]
where
\[
S(r,cn;q)=\sum_{a\bmod q}^{*}e((ra+cn\overline a)/q).
\]
This is the classical Kloosterman sum. Next using Weil's bound
(see \cite[Chapter~11]{IwaniecKowalski2004}),\eqref{hm:eq:deligne} and taking sums over \(n\) gives
\begin{equation}\label{hm:eq-poisson-block}
E(N,M;c)
\ll_{f,B}\left(q\sqrt{\frac NM}+\sqrt{\frac{NM}{q}}\right)
\log(2N)(1+P/Q)^{-B}.
\end{equation}

Finally, combining \eqref{hm:eq-wilton-block} 
\eqref{hm:eq-poisson-block}, and using
$\min(x,y)\le x^{2/3}y^{1/3}$, yields
\begin{equation}\label{hm:eq-interpolated-block}
\begin{aligned}
E(N,M;c)
&\ll_{f,B}
\left[\min\!\left(q^{1/2}M^{1/2},q\sqrt{N/M}\right)
+\sqrt{NM/q}\right]\log(2P)(1+P/Q)^{-B}\\
&\ll_{f,B}\left[q^{2/3}P^{1/6}+\sqrt{P/q}\right]
\log(2P)(1+P/Q)^{-B}.
\end{aligned}
\end{equation}
Summing over the dyadic blocks gives
\begin{equation}\label{hm:eq-dual-error}
R_{1q}(h, l)
\ll_{f,\varepsilon}
q^\varepsilon\left(q^{2/3}Q^{1/6}+q^{-1/2}Q^{1/2}\right)
\ll_{f,\varepsilon}q^{11/12+\varepsilon}.
\end{equation}
Combining \eqref{hm:eq-diagonal-main}, \eqref{hm:eq-direct-error},
\eqref{hm:eq-dual-correction}, and \eqref{hm:eq-dual-error}, this completes the proof of 
\eqref{hm:eq-weighted-first}.
\end{proof}

\section{Outline of the proof}
\label{hm:sec:outline}

Fix a positive integer $k$, and set
\[
 X=q^{1/(32k)},\qquad \mathcal L=\log X.
\]
For $u>0$, define
\[
 V(u)=\mathbf 1_{u\le X}\frac{\log(X/u)}{\log X},
\]
and put
\begin{equation}\label{hm:eq:A}
 A(\chi)=\sum_{n,m\ge1}
 \frac{\lambda_f(n)\chi(n)\overline{\chi(m)}}{\sqrt{nm}}V(nm).
\end{equation}
Define
\begin{equation}\label{hm:eq:S1}
 S_1:=\sum_{\chi\in\mathcal X_q^+}
 F(\chi)A(\chi)^{k-1}\overline{A(\chi)}^{\,k},
\end{equation}
and
\begin{equation}\label{hm:eq:S2}
 S_2:=\sum_{\chi\in\mathcal X_q^+}|A(\chi)|^{2k}.
\end{equation}
Since $L(1/2,\overline\chi)=\overline{L(1/2,\chi)}$, H\"older's inequality
and the inclusion of $\mathcal X_q^+$ in the primitive-character family give
\begin{equation}\label{hm:eq:moment-reduction}
 \sum_{\chi\bmod q}^{*}
 \left|L\!\left(\tfrac12,f\otimes\chi\right)
 L\!\left(\tfrac12,\chi\right)\right|^{2k}
 \ge \frac{|S_1|^{2k}}{S_2^{2k-1}}.
\end{equation}

The estimates to be established are
\begin{equation}\label{hm:eq:S1-size}
 S_1\asymp_{f,k}q(\log q)^{2k^2}
\end{equation}
and
\begin{equation}\label{hm:eq:S2-size}
 S_2\asymp_{f,k}q(\log q)^{2k^2}.
\end{equation}

Together with \eqref{hm:eq:moment-reduction}, these estimates yield
the claimed lower bound for every fixed positive integer $k$.

\section{Estimates for \texorpdfstring{$S_2$}{S2}}
\label{hm:sec:S2}

Expanding $|A(\chi)|^{2k}$ gives
\[
 S_2=\sum_{\chi\in\mathcal X_q^+}
 \sum_{\boldsymbol n,\boldsymbol m,\boldsymbol\ell,\boldsymbol r\ge1}
 \left(\prod_{i=1}^k
 \frac{\lambda_f(n_i)\lambda_f(m_i)}{\sqrt{n_im_i\ell_ir_i}}
 V(n_i\ell_i)V(m_ir_i)\right)
 \prod_{i=1}^k
 \chi(n_i)\chi(r_i)\overline{\chi(\ell_i)\chi(m_i)}.
\]
Here and below, a sum over vectors of positive integers means that all
their components are summed independently. Put
\[
 N=\prod_{i=1}^k n_ir_i,\qquad M=\prod_{i=1}^k \ell_im_i.
\]
On the support of the weights,
\[
 NM=\prod_{i=1}^k(n_i\ell_i)(m_ir_i)\le X^{2k}=q^{1/16}.
\]
In particular, $N,M<q/2$ for sufficiently large $q$, and all variables
are coprime to $q$. Thus $N\equiv M\pmod q$ is equivalent to $N=M$,
while $N\equiv-M\pmod q$ has no solutions. Then using character's orthogonality yields
\begin{equation}\label{hm:eq:S2-diagonal}
 S_2=\frac{q-1}{2}M_{21}+S_{22},
 \qquad S_{22}:=-|A(\chi_0)|^{2k}.
\end{equation}
Here
\begin{equation}\label{hm:eq:M21}
 M_{21}:=
 \sum_{\substack{\boldsymbol n,\boldsymbol m,\boldsymbol\ell,\boldsymbol r\ge1\\N=M}}\prod_{i=1}^k
 \frac{\lambda_f(n_i)\lambda_f(m_i)}{\sqrt{n_im_i\ell_ir_i}}
 V(n_i\ell_i)V(m_ir_i).
\end{equation}
By \eqref{hm:eq:deligne}, we obtain
\begin{equation}\label{hm:eq:S22}
 \begin{aligned}
 |S_{22}|=|A(\chi_0)|^{2k}\le\left(\sum_{nm\le X}\frac{d(n)}{\sqrt{nm}}\right)^{2k}
 =\left(\sum_{u\le X}\frac{d_3(u)}{\sqrt u}\right)^{2k}\ll_{k,\varepsilon}X^{k+\varepsilon}
 \ll_{k,\varepsilon}q^{1/32+\varepsilon}.
 \end{aligned}
\end{equation}

The Mellin inversion gives
\begin{equation}\label{hm:eq:V-mellin}
 V(u)=\frac1{2\pi i\mathcal L}
 \int_{(c)}\left(\frac Xu\right)^s\frac{ds}{s^2},
 \qquad u>0,\quad c>0,
\end{equation}
where $\int_{(c)}$ denotes integration upwards along $\Re s=c$.
Define the multiple Dirichlet series
\[
 \mathcal F_k(\boldsymbol s,\boldsymbol t)
 :=\sum_{\substack{\boldsymbol n,\boldsymbol m,\boldsymbol\ell,\boldsymbol r\ge1\\N=M}}\prod_{i=1}^k
 \frac{\lambda_f(n_i)\lambda_f(m_i)}
 {(n_i\ell_i)^{1/2+s_i}(m_ir_i)^{1/2+t_i}}.
\]
This series converges absolutely when all $\Re s_i,\Re t_i>0$.
Applying \eqref{hm:eq:V-mellin} to all
$2k$ weights and interchanging the absolutely convergent sums and
integrals, we obtain
\begin{equation}\label{hm:eq:M21-mellin}
 M_{21}=\frac1{(2\pi i)^{2k}\mathcal L^{2k}}
 \int_{(c)^{2k}}
 \mathcal F_k(\boldsymbol s,\boldsymbol t)
 X^{\sum_{i=1}^k(s_i+t_i)}
 \prod_{i=1}^k\frac{ds_i\,dt_i}{s_i^2t_i^2},
 \qquad c=\mathcal L^{-1},
\end{equation}
where $\int_{(c)^{2k}}$ denotes the product of $2k$ vertical
integrals.

After calculation of Euler-product gives
\begin{equation}\label{hm:eq:F-factorization}
 \mathcal F_k(\boldsymbol s,\boldsymbol t)
 =H(\boldsymbol s,\boldsymbol t)G(\boldsymbol s,\boldsymbol t)
 \prod_{i,j=1}^k P_f(1+s_i+t_j),
 \qquad P_f(s):=\zeta(s)L(s,f\times f).
\end{equation}
For $0<\eta<1/4$ and $\Re s_i,\Re t_i\ge-\eta$, where
\begin{equation}\label{hm:eq:H-product}
 \begin{aligned}
 H(\boldsymbol s,\boldsymbol t)
 :=\prod_p
 \frac{\mathcal F_{k,p}(\boldsymbol s,\boldsymbol t)}
 {G_p(\boldsymbol s,\boldsymbol t)
  \prod_{i,j=1}^kP_{f,p}(1+s_i+t_j)}
 =\prod_p\bigl(1+O_{k,\eta}(p^{-2+4\eta})\bigr),
 \end{aligned}
\end{equation}
and
\begin{equation}\label{hm:eq:G-regular}
 G(\boldsymbol s,\boldsymbol t)
 :=\prod_{i,j=1}^k
 L(1+s_i+s_j,f)L(1+t_i+t_j,f).
\end{equation}

To deal with \eqref{hm:eq:M21-mellin}, substitute
$s_i=z_i/\mathcal L$ and $t_i=w_i/\mathcal L$. Hence,
\begin{equation}\label{hm:eq:M21-rescaled}
 \begin{aligned}
 M_{21}=\frac1{(2\pi i)^{2k}}
 \int_{(1)^{2k}}\mathcal F_k\!\left(\frac{\boldsymbol z}{\mathcal L},
                     \frac{\boldsymbol w}{\mathcal L}\right)
 e^{\sum_i^k(z_i+w_i)}\prod_{i=1}^k\frac{dz_i\,dw_i}{z_i^2w_i^2}.
 \end{aligned}
\end{equation}

Put $D=2k^2$ and $T=\sqrt{\mathcal L}$. Truncate the integral at
$|\Im z_i|,|\Im w_i|\le T$ before expanding the integrand.
The identity
\[
 L(s,f\times f)=\zeta(2s)\sum_{n\ge1}\frac{\lambda_f(n)^2}{n^s}
 \qquad(\Re s>1)
\]
and its nonnegative coefficients imply
$|L(1+2/\mathcal L+it,f\times f)|\ll_f\mathcal L$.
Similarly, $|\zeta(1+2/\mathcal L+it)|\ll\mathcal L$.

Then for any $\varepsilon>0$  in terms of $k$ and using
\eqref{hm:eq:F-factorization}, we obtain
\[
 \left|\mathcal F_k\!\left(
 \frac{\boldsymbol z}{\mathcal L},
 \frac{\boldsymbol w}{\mathcal L}\right)\right|
 \ll_{f,k}\mathcal L^D
 \prod_{i=1}^k
 \left(1+\frac{|\Im z_i|}{\mathcal L}\right)^{1/4}
 \left(1+\frac{|\Im w_i|}{\mathcal L}\right)^{1/4}
.
\]
After multiplication by the integration, the contribution from
$\max_i(|\Im z_i|,|\Im w_i|)>T$ is
$O_{f,k}(\mathcal L^D/T)=O_{f,k}(\mathcal L^{D-1/2})$.

Writing $\int_{(1;T)^{2k}}$ for the $2k$ times product of integral where the segments from
$1-iT$ to $1+iT$, we conclude that
\begin{equation}\label{hm:eq:M21-truncated}
 \begin{aligned}
 M_{21}=\frac1{(2\pi i)^{2k}}
 \int_{(1;T)^{2k}}\mathcal F_k\!\left(\frac{\boldsymbol z}{\mathcal L},
                     \frac{\boldsymbol w}{\mathcal L}\right)
 e^{\sum_i^k(z_i+w_i)}\prod_{i=1}^k\frac{dz_i\,dw_i}{z_i^2w_i^2}+O_{f,k}(\mathcal L^{D-1/2}).
 \end{aligned}
\end{equation}

Expanding the pole factors on the truncated contours and extending the
limiting kernel to the full contours gives
\begin{equation}\label{eq:M21-leading-integral}
 M_{21}=B_{f,k}C_k\mathcal L^{2k^2}
 +O_{f,k}(\mathcal L^{2k^2-1/2}),
\end{equation}
where
\[
R_f=\operatorname*{Res}_{s=1}L(s,f\times f)>0, \qquad
B_{f,k}:=H(\mathbf{0},\mathbf{0})G(\mathbf{0},\mathbf{0})R_f^{k^2}.
\]
\[
C_k:=\frac{1}{(2\pi i)^{2k}} \int_{(1)^{2k}} \frac{e^{\sum_i^k(z_i+w_i)}}{\prod_{i,j=1}^k(z_i+w_j)^2} \prod_{i=1}^k\frac{dz_i\,dw_i}{z_i^2w_i^2}.
\]

By the similar argument for the constant $C$ in \cite[p.~12]{HeapLiZhao2022}, one has 
$C_k, B_{f,k}>0$. Combining \eqref{hm:eq:S2-diagonal} and \eqref{hm:eq:S22} yields
\[
 S_2=\frac{q-1}{2}B_{f,k}C_k(\log X)^{2k^2}
 +O_{f,k,\varepsilon}\!\left(
 q(\log X)^{2k^2-1/2}+q^{1/32+\varepsilon}\right).
\]

Therefore, we obtain
$S_2\asymp_{f,k}q(\log q)^{2k^2}$. This completes the proof of 
\eqref{hm:eq:S2-size}.

\section{Estimates for \texorpdfstring{$S_1$}{S1}}
\label{hm:sec:S1}

To expand the auxiliary polynomial, use pairs $(a_i,b_i)$,
$1\le i\le k-1$, from $A(\chi)^{k-1}$ and pairs $(c_j,d_j)$,
$1\le j\le k$, from $\overline{A(\chi)}^{\,k}$.
For convinence, we denote
\[
 \begin{gathered}
 h=\prod_{i=1}^{k-1}a_i\prod_{j=1}^k d_j,
 \qquad
 \ell=\prod_{i=1}^{k-1}b_i\prod_{j=1}^k c_j,
 \qquad g=(h,\ell),\\
 \Lambda=\prod_{i=1}^{k-1}\lambda_f(a_i)
          \prod_{j=1}^k\lambda_f(c_j),
 \qquad
 \mathcal V=\prod_{i=1}^{k-1}V(a_ib_i)
            \prod_{j=1}^kV(c_jd_j).
 \end{gathered}
\]
\begin{equation}\label{hm:eq:S1-tuple}
 \begin{aligned}
 S_1&=\sum_{1\le\boldsymbol a,\boldsymbol b,\boldsymbol c,\boldsymbol d\le X}
 \frac{\Lambda\mathcal V}{\sqrt{h\ell}}\,T_q(h,\ell),\\
 T_q(h,\ell)&:=\sum_{\chi\in\mathcal X_q^+}
 F(\chi)\chi(h)\overline{\chi(\ell)}.
 \end{aligned}
\end{equation}
Here and below, a sum over vectors of positive integers means that all
their components are summed independently. On the support of $\mathcal V$,
$h\ell\le X^{2k-1}<q^{1/16}$.
The absolute sum of the short coefficients satisfies
\begin{equation}\label{hm:eq:S1-coefficient-norm}
 \begin{aligned}
 \sum_{1\le\boldsymbol a,\boldsymbol b,\boldsymbol c,\boldsymbol d\le X}
 \frac{|\Lambda|\mathcal V}{\sqrt{h\ell}}=\left(\sum_{nm\le X}
        \frac{|\lambda_f(n)|V(nm)}{\sqrt{nm}}\right)^{2k-1}\le\left(\sum_{u\le X}\frac{d_3(u)}{\sqrt u}\right)^{2k-1}
 \ll_{k,\varepsilon}X^{k-1/2+\varepsilon}.
 \end{aligned}
\end{equation}
By Lemma~\ref{hm:lem:weighted}, we obtain
\begin{equation}\label{hm:eq:S1-main-reduction}
 S_1=\frac{q-1}{2}E_q(1)M_k(X)
 +O_{f,k,\varepsilon}(q^{\theta_k+\varepsilon}),
\end{equation}
where $\theta_k=91/96-1/(64k)<1$ and
\begin{align}
E_q(s) &= 1-\frac{\lambda_f(q)}{q^s}+\frac{1}{q^{2s}}, \label{eq:Eq-definition} \\
M_k(X) &:= L(1,f)\sum_{1\le\mathbf{a},\mathbf{b},\mathbf{c},\mathbf{d}\le X} \Lambda\mathcal{V}\,\frac{g}{h\ell}C_{\ell/g}(1), \label{eq:Mk-definition} \\
C_b(s) &= \sum_{d\mid b}\frac{\mu(d)\lambda_f(b/d)}{d^s}. \label{eq:Cb-definition}
\end{align}

To evaluate this signed main term, introduce $k-1$ Mellin variables
$\boldsymbol s=(s_1,\ldots,s_{k-1})$ and $k$ variables
$\boldsymbol t=(t_1,\ldots,t_k)$. Define
\begin{equation}\label{hm:eq:R-series}
 \begin{aligned}
 \mathcal R_k(\boldsymbol s,\boldsymbol t)
 :=\sum_{1\le\boldsymbol a,\boldsymbol b,\boldsymbol c,\boldsymbol d\le X}
 \Lambda\frac{g}{h\ell}C_{\ell/g}(1)
 \prod_{i=1}^{k-1}(a_ib_i)^{-s_i}
 \prod_{j=1}^k(c_jd_j)^{-t_j}.
 \end{aligned}
\end{equation}
Mellin inversion and absolute convergence give
\begin{equation}\label{hm:eq:Mk-mellin}
 \begin{aligned}
 M_k(X)=\frac{L(1,f)}{(2\pi i)^{2k-1}\mathcal L^{2k-1}}
 \int_{(c)^{2k-1}}
 \mathcal R_k(\boldsymbol s,\boldsymbol t)
 X^{\sum_{i=1}^{k-1} s_i+\sum_{j=1}^{k}t_j}\prod_{i=1}^{k-1}\frac{ds_i}{s_i^2}
 \prod_{j=1}^k\frac{dt_j}{t_j^2},
 \qquad c=\mathcal L^{-1}.
 \end{aligned}
\end{equation}

Write $\lambda=\lambda_f(p)$ and fix $0<\eta<1/4$. Expanding the local Euler factors gives
\begin{equation}\label{hm:eq:R-local-expansion}
 \begin{aligned}
 \mathcal R_{k,p}(\boldsymbol s,\boldsymbol t)
 ={}&1+(\lambda^2+1)\sum_{j=1}^kp^{-1-t_j}
      +2\lambda\sum_{i=1}^{k-1}p^{-1-s_i}\\
 &+(\lambda^2+1)\sum_{i=1}^{k-1}\sum_{j=1}^k
       p^{-1-s_i-t_j}\\
 &+\lambda\sum_{i,a=1}^{k-1}p^{-1-s_i-s_a}
  +\lambda\sum_{j,b=1}^kp^{-1-t_j-t_b}
  +O_{k,\eta}(p^{-2+4\eta}),
 \end{aligned}
\end{equation}
uniformly for $\Re s_i,\Re t_j\ge-\eta$.

Comparing the
coefficients in \eqref{hm:eq:R-local-expansion}, we obtain
\begin{equation}\label{hm:eq:R-factorization}
 \begin{aligned}
 \mathcal R_k(\boldsymbol s,\boldsymbol t)
 =H_1(\boldsymbol s,\boldsymbol t)
 K_f(\boldsymbol s,\boldsymbol t)
 \prod_{j=1}^kP_f(1+t_j)\prod_{i=1}^{k-1}\prod_{j=1}^kP_f(1+s_i+t_j),
 \end{aligned}
\end{equation}
where
\[
 K_f(\boldsymbol s,\boldsymbol t)
 =\prod_{i=1}^{k-1}L(1+s_i,f)^2
  \prod_{i,a=1}^{k-1}L(1+s_i+s_a,f)
  \prod_{j,b=1}^kL(1+t_j+t_b,f).
\]

To deal with \eqref{hm:eq:Mk-mellin}, make the substitutions
$s_i=z_i/\mathcal L$ and $t_j=w_j/\mathcal L$.
An argument analogous to that for $S_2$ gives
\begin{equation}\label{hm:eq:Mk-asymptotic}
 M_k(X)=\mathfrak a_{f,k}J_k(\log X)^{2k^2}
 +o_{f,k}\bigl((\log X)^{2k^2}\bigr),
\end{equation}
where
\[
 \mathfrak a_{f,k}
 :=H_1(\boldsymbol0,\boldsymbol0)R_f^{k^2}L(1,f)^{2k^2},
\]
and
\begin{equation}\label{hm:eq:Jk-complex}
 J_k:=\frac1{(2\pi i)^{2k-1}}
 \int_{(1)^{2k-1}}
 \frac{e^{\sum_{i=1}^{k-1} z_i+\sum_{j=1}^{k}w_j}}
 {\prod_{i=1}^{k-1}\prod_{j=1}^k(z_i+w_j)^2}
 \prod_{i=1}^{k-1}\frac{dz_i}{z_i^2}
 \prod_{j=1}^k\frac{dw_j}{w_j^4}.
\end{equation}

The same positivity argument as for $C_k$ gives $J_k>0, \mathfrak a_{f,k}>0$. Combining \eqref{hm:eq:S1-main-reduction} and
\eqref{hm:eq:Mk-asymptotic} gives
\begin{equation}\label{hm:eq:S1-asymptotic}
 S_1=\frac{q-1}{2}\mathfrak a_{f,k}J_k(\log X)^{2k^2}
 +o_{f,k}\bigl(q(\log X)^{2k^2}\bigr).
\end{equation}
This proves \eqref{hm:eq:S1-size}. 

\section*{Declaration of AI use}

During the preparation of this work, the authors used OpenAI's
GPT-6 Astra to assist with calculations, literature searches, and
typesetting. The authors reviewed and edited the resulting material
and take responsibility for the accuracy, originality, and integrity
of the manuscript.

\end{document}